\documentclass{article}
\usepackage{amsmath}
\usepackage{amssymb}
\usepackage{amscd}

\def\be{\begin{equation}}       \def\ee{\end{equation}}
\def\bd{\begin{displaymath}}    \def\ed{\end{displaymath}}
\def\beq{\begin{eqnarray}}      \def\eeq{\end{eqnarray}}
\def\bseq{\begin{eqnarray*}}    \def\eseq{\end{eqnarray*}}
\def\ba{\begin{array}}          \def\ea{\end{array}}
\def\ben{\begin{enumerate}}     \def\een{\end{enumerate}}

\def\cop{\Delta}

\def\cnt{\varepsilon}
\def\ot{\otimes}

\def\GLqtwo{{GL}_{q}(2)}

\def\GLpqtwo{{GL}_{p,q}(2)}
\def\GLhh'two{{GL}_{h,h'}(2)}

\def\Grss'{{G}_{r}^{s,s'}}
\def\Gmkk'{{G}_{m}^{k,k'}}

\def\GLtwo{{GL}({2})}

\def\ident{{\bf 1}}

\def\ca{{\mathcal A}}

\def\crr{{\mathcal R}}

\def\ms{\mathsf s}

\numberwithin{equation}{section}
\begin{document}
\begin{center}
{\large \bf 
DUALITY FOR COLOURED QUANTUM GROUPS \footnote{Max-Planck-Institute
Preprint MPI-MIS-64/2000, Edinburgh University Preprint MS-00-006}
}

\bigskip
\bigskip
\bigskip
\bigskip
\bigskip
\bigskip

{\large\bf Deepak Parashar
}

\smallskip
{\it Max-Planck-Institute for Mathematics in the Sciences\\
Inselstrasse 22-26, D-04103 Leipzig\\ 
Germany\\
}

\smallskip
{\tt 
e-mail: Deepak.Parashar@mis.mpg.de
}

\bigskip
\bigskip
\bigskip
\bigskip
{\large \bf Abstract}
\end{center}
\medskip

Duality between the coloured quantum group and the coloured quantum
algebra corresponding to $GL(2)$ is established. The coloured
$L^{\pm}$ functionals are constructed and the dual algebra is derived
explicitly. These functionals are then employed to give a coloured
generalisation of the differential calculus on quantum $GL(2)$ within the
framework of the $R$-matrix approach.
\par

\bigskip
\bigskip
\bigskip
\begin{center}
{\it Lett. Math. Phys.} {\bf 53} (2000) 29-40
\end{center}

\newpage

\section{Introduction}

It is well-known [1-3] that the standard quantum group relations can be
extended
by parametrising the corresponding generators using some continuous
`colour' variables and redefining the associated algebra and coalgebra in
a way that all Hopf algebraic properties remain preserved. This results 
in a `coloured' extension of a quantum group and the associated $R$-matrix
is a nonadditive type $R$-matrix.
\\
Non-additive type solutions of the Yang-Baxter equation were
first discovered in the study of integrable models by Bazhanov and
Stroganov [4]. In the context of quantum groups, a coloured $R$-matrix
solution was first introduced by Murakami [5] and subsequently in [6,7].
In the context of knot theory, Ohtsuki [8] introduced coloured
quasitriangular Hopf algebras characterised by the existence of a coloured
universal  $R$-matrix. Jordanian deformations also admit coloured
extensions [3,9,10] and the associated $R$-matrix is `colour' triangular
{\it i.e.}, a coloured extension of the notion of triangularity. More
recently, it has been shown [10] that coloured Jordanian deformations can
be obtained from the coloured quantum deformations by means of a
contraction procedure. Coloured generalisations of quantum groups can also
be understood as an application of the twisting procedure, in a manner
similar to the multiparameter generalisation of quantum groups.
\\
Even though a considerable interest has been generated in the coloured
generalisation of quantum groups, some aspects of their basic algebraic
and geometric structure still need thorough investigation. As is the case
with ordinary quantum groups, one works on the coloured
generalisation of {\em either} the quantised algebra of functions on a
group {\em or} that of the quantised universal enveloping algebra. It is 
already well-established [11,12] that the two structures are dual to each
other in the uncoloured case. The purpose of this letter is {\em
two-fold}. {\em Firstly}, we address
the problem of duality between coloured quantum groups and coloured
quantum algebras. This would then provide us with some more information
about the algebraic structure underlying these objects. {\em Secondly}, it
is well-known that in the uncoloured case the problem of duality
naturally leads us to
formulate a differential calculus on quantum groups. In the same spirit,
we investigate the geometric structure of coloured quantum groups by
the construction of a differential calculus. Throughout this letter, we
shall focus on the coloured generalisation of the simplest
single-parameter quantum deformation of $GL(2)$, denoted as
$GL_{q}^{\lambda,\mu}(2)$. In pursuing our aim, we shall adhere to the
convenient $R$-matrix approach [13,14]. 
\par
In Section $2$, we recall the definition of the coloured quantum group
$GL^{\lambda , \mu}_{q}(2)$. In Section $3$, we construct the coloured
$L^{\pm}$ functionals for $GL_{q}^{\lambda,\mu}(2)$ and derive explicitly
its dual algebra {\it i.e.} the coloured quantised universal enveloping
algebra and also exhibit its Hopf algebraic structure. Section $4$ is
devoted to a coloured generalisation of the $R$-matrix approach (also
known as the {\em constructive} procedure) to construct a differential
calculus on the coloured quantum group $GL_{q}^{\lambda,\mu}(2)$. The
Letter concludes with a summary of the results in Section $5$.

\section{The Coloured Quantum Group $GL_{q}^{\lambda,\mu}(2)$}

For the case of a single-parameter quantum deformation of $\GLtwo$ (with
deformation parameter $q$), its `coloured' version [1] is given by the
$R$-matrix,
\begin{equation}
R_{q}^{\lambda,\mu} = \begin{pmatrix}
q^{1-(\lambda-\mu)} & 0 & 0 & 0\\
0 & q^{\lambda+\mu} & 0 & 0\\
0 & q-q^{-1} & q^{-(\lambda+\mu)} & 0\\
0 & 0 & 0 & q^{1+(\lambda-\mu)}\
\end{pmatrix}
\end{equation}

which satisfies
\begin{equation}
R_{12}^{\lambda,\mu}R_{13}^{\lambda,\nu}R_{23}^{\mu,\nu} =
R_{23}^{\mu,\nu}R_{13}^{\lambda,\nu}R_{12}^{\lambda,\mu}
\end{equation}
the so-called `coloured' quantum Yang-Baxter equation (CQYBE). The
coloured $R$-matrix provides a
nonadditive-type solution $R^{\lambda,\mu} \neq R(\lambda - \mu)$ of the
Yang-Baxter equation, which is in general
multicomponent and the parameters $\lambda$, $\mu$, $\nu$ are considered
as `colour' parameters. This gives rise to the coloured $RTT$ relations 
\begin{equation}
R_{q}^{\lambda,\mu}T_{1\lambda}T_{2\mu}=T_{2\mu}T_{1\lambda}
R_{q}^{\lambda,\mu}
\end{equation}
(where $T_{1\lambda}=T_{\lambda}\ot \ident$ and
$T_{2\mu}=\ident\ot T_{\mu}$) in which the entries of the $T$
matrices carry colour  dependence {\sl i.e.}
$T_{\lambda}= \left(
\begin{smallmatrix}a_{\lambda}&b_{\lambda}\\c_{\lambda}&d_{\lambda}
\end{smallmatrix}\right)$,
$T_{\mu}= \left(
\begin{smallmatrix}a_{\mu}&b_{\mu}\\c_{\mu}&d_{\mu}\end{smallmatrix} 
\right)$.
The coproduct and counit for the coalgebra structure are given by
$\cop (T_{\lambda})=T_{\lambda}\dot{\ot} T_{\lambda}$,
$\cnt (T_{\lambda})=\ident$
and depend on one colour parameter at a time. By contrast, the algebra 
structure is more complicated with generators of two different 
colours appearing simultaneously in the algebraic relations. The full Hopf
algebraic structure can be constructed resulting in a coloured extension
of the quantum group. Since $\lambda$ and $\mu$ are continuous variables,
this implies the coloured quantum group has an infinite number of
generators. The quantum determinant
$D_{\lambda}=a_{\lambda}d_{\lambda}-r^{-(1+2\lambda)}c_{\lambda}b_{\lambda}$
is group-like but not central, even in the case of the single-parameter
coloured deformation. This is a crucial difference with the ordinary
uncoloured quantum groups. The antipode is given as
\begin{equation}
S(T_{\lambda})= 
D_{\lambda}^{-1}\begin{pmatrix}d_{\lambda}&-r^{1+2\lambda}b_{\lambda}\\
-r^{-1-2\lambda}c_{\lambda}&a_{\lambda}\end{pmatrix}
\end{equation}
This results in a coloured generalisation of a quantum group in the
framework of the $FRT$ formalism.
It should be noted that while the colourless limit $\lambda = \mu = 0$
gives back the ordinary single-parameter deformed quantum group, the
monochromatic limit $\lambda = \mu \neq 0$ gives rise to the uncoloured
two-parameter deformed quantum group. Analogous to the quantum group
case, one can also define the coloured braid group representation as
\begin{equation}
\hat{\crr}^{\lambda , \mu} = P \crr^{\lambda , \mu}
\end{equation}
where $P$ is the usual permutation matrix. This coloured braid group
representation turns out to be a solution of the coloured braided
Yang-Baxter equation
\begin{equation}
\hat{R}_{23}^{\lambda,\mu} \hat{R}_{12}^{\lambda,\nu}
\hat{R}_{23}^{\mu,\nu} =
\hat{R}_{12}^{\mu,\nu} \hat{R}_{23}^{\lambda,\nu}
\hat{R}_{12}^{\lambda,\mu}
\end{equation}

\section{The Dual Algebra for $GL_{q}^{\lambda,\mu}(2)$}

Here we derive explicitly the algebra dual to the coloured
quantum group $GL_{q}^{\lambda,\mu}(2)$. If we denote generators of the
yet unknown dual algebra by 
$\{ A_{\lambda},B_{\lambda},C_{\lambda},D_{\lambda} \}$ and
$\{ A_{\mu},B_{\mu},C_{\mu},D_{\mu} \}$, then the following
pairings hold
\begin{eqnarray}
&&\langle A_{\lambda}, T_{\lambda} \rangle = 
\langle A_{\lambda}, T_{\mu} \rangle =
\langle A_{\mu}, T_{\lambda} \rangle =
\langle A_{\mu}, T_{\mu} \rangle = 
\left( \begin{smallmatrix}1&0\\0&0\end{smallmatrix} \right)\\
&&\langle B_{\lambda}, T_{\lambda} \rangle =
\langle B_{\lambda}, T_{\mu} \rangle = 
\langle B_{\mu}, T_{\lambda} \rangle =
\langle B_{\mu}, T_{\mu} \rangle =
\left( \begin{smallmatrix}0&1\\0&0\end{smallmatrix} \right)\\
&&\langle C_{\lambda}, T_{\lambda} \rangle =
\langle C_{\lambda}, T_{\mu} \rangle = 
\langle C_{\mu}, T_{\lambda} \rangle =
\langle C_{\mu}, T_{\mu} \rangle =
\left( \begin{smallmatrix}0&0\\1&0\end{smallmatrix} \right)\\
&&\langle D_{\lambda}, T_{\lambda} \rangle =
\langle D_{\lambda}, T_{\mu} \rangle =
\langle D_{\mu}, T_{\lambda} \rangle =
\langle D_{\mu}, T_{\mu} \rangle =
\left( \begin{smallmatrix}0&0\\0&1\end{smallmatrix} \right)
\end{eqnarray}
where $T_{\lambda}$ and $T_{\mu}$ are the matrices of generators of
$GL_{q}^{\lambda,\mu}(2)$.

\subsection{Coloured $L^{\pm}$ functionals}

Recall [11,14] that the $L^{\pm}$ functionals for the uncoloured
single-parameter quantum group $\GLqtwo$ are expressed in terms of the
matrices
\begin{equation}
L^{+} = c^{+} q^{1/2} \begin{pmatrix}
q^{H/2} & q^{-1/2}(q-q^{-1})X^{+}\\
0 & q^{-H/2}
\end{pmatrix}
\end{equation}
\begin{equation}
L^{-} = c^{-} q^{-1/2} \begin{pmatrix}
q^{-H/2} & 0\\
q^{1/2}(q^{-1}-q)X^{-} & q^{H/2}
\end{pmatrix} 
\end{equation}
where $c^{+}$ and $c^{-}$ are free parameters and $H=A-D$, $X^{+}=C$,
$X^{-}=B$ are the generators of the algebra dual to $\GLqtwo$. The
$R$-matrix for the coloured quantum group $GL_{q}^{\lambda,\mu}(2)$
can be written as
\begin{equation}
R_{12}
= q^{1/2} \begin{pmatrix}
q^{-1/2}q^{1-\lambda+\mu} & 0 & 0 & 0\\
0 & q^{-1/2}q^{\lambda+\mu} & 0 & 0\\ 
0 & q^{-1/2}(q-q^{-1}) & q^{-1/2}q^{-(\lambda+\mu)} & 0\\
0 & 0 & 0 & q^{-1/2}q^{1+\lambda-\mu}
\end{pmatrix}
\end{equation}
The corresponding $R^{+}$ and $R^{-}$ matrices read
\begin{eqnarray}  
R^{+} &=& c^{+}R_{21} \notag \\
&=& c^{+} q^{1/2} \begin{pmatrix}
q^{-1/2}q^{1-\lambda+\mu} & 0 & 0 & 0\\
0 & q^{-1/2}q^{-(\lambda+\mu)} & q^{-1/2}(q-q^{-1}) & 0\\
0 & 0 & q^{-1/2}q^{\lambda+\mu} & 0\\
0 & 0 & 0 & q^{-1/2}q^{1+\lambda-\mu} 
\end{pmatrix}
\end{eqnarray}
\begin{eqnarray}
R^{-} &=& c^{-}R_{12}^{-1} \notag \\
&=& c^{-} q^{-1/2} \begin{pmatrix}
q^{1/2}q^{-(1-\lambda+\mu)} & 0 & 0 & 0\\
0 & q^{1/2}q^{-(\lambda+\mu)} & 0 & 0\\  
0 & -q^{1/2}(q-q^{-1}) & q^{1/2}q^{\lambda+\mu} & 0\\
0 & 0 & 0 & q^{1/2}q^{-(1+\lambda-\mu)}
\end{pmatrix}
\end{eqnarray} 
The coloured $L^{\pm}$ functionals can be expressed as
\begin{eqnarray}
L^{+}_{\lambda} &=& c^{+} q^{1/2} \begin{pmatrix}
q^{H_{\lambda}/2}q^{\mu H_{\lambda}-\lambda H'_{\lambda}} &
q^{-1/2}(q-q^{-1})C_{\lambda}\\
0 & q^{-H_{\lambda}/2}q^{\mu H_{\lambda}+\lambda H'_{\lambda}}
\end{pmatrix}\\
L^{+}_{\mu} &=& c^{+} q^{1/2} \begin{pmatrix}
q^{H_{\mu}/2}q^{\mu H_{\mu}-\lambda H'_{\mu}} & 
q^{-1/2}(q-q^{-1})C_{\mu}\\
0 & q^{-H_{\mu}/2}q^{\mu H_{\mu}+\lambda H'_{\mu}}
\end{pmatrix}
\end{eqnarray}
\begin{eqnarray}
L^{-}_{\lambda} &=& c^{-} q^{-1/2} \begin{pmatrix}
q^{-H_{\lambda}/2}q^{\lambda H_{\lambda}-\mu H'_{\lambda}} & 0\\
q^{1/2}(q^{-1}-q)B_{\lambda} & 
q^{H_{\lambda}/2}q^{\lambda H_{\lambda}+\mu H'_{\lambda}}
\end{pmatrix}\\
L^{-}_{\mu} &=& c^{-} q^{-1/2} \begin{pmatrix}
q^{-H_{\mu}/2}q^{\lambda H_{\mu}-\mu H'_{\mu}} & 0\\
q^{1/2}(q^{-1}-q)B_{\mu} & 
q^{H_{\mu}/2}q^{\lambda H_{\mu}+\mu H'_{\mu}}
\end{pmatrix}
\end{eqnarray}
where $H_{\lambda}=A_{\lambda}-D_{\lambda}$,
$H'_{\lambda}=A_{\lambda}+D_{\lambda}$ and
$H_{\mu}=A_{\mu}-D_{\mu}$, $H'_{\mu}=A_{\mu}+D_{\mu}$.
Note that each one of $L^{\pm}_{\lambda}$ and $L^{\pm}_{\mu}$ depends on
{\em both} the
colour parameters $\lambda$ and $\mu$. The notation $L^{\pm}_{\lambda}$
implies that the generators of the dual carry $\lambda$ dependence, and
similarly $L^{\pm}_{\mu}$ implies that the generators of the dual carry
$\mu$ dependence. The duality pairings are then given by the action of
the functionals $L^{\pm}_{\lambda}$ and $L^{\pm}_{\mu}$ on the
$T$-matrices $T_{\lambda}$ and $T_{\mu}$ of the coloured quantum group
$GL_{q}^{\lambda,\mu}(2)$
\begin{eqnarray}
&&(L^{+}_{\lambda|\mu})^{a}_{b} (T_{\lambda|\mu})^{c}_{d} =
(R^{+})^{ac}_{bd}\\
&&(L^{-}_{\lambda|\mu})^{a}_{b} (T_{\lambda|\mu})^{c}_{d} =
(R^{-})^{ac}_{bd}
\end{eqnarray}
The notation $\lambda|\mu$ in the subscript in the above relations means
{\em either} $\lambda$ {\em or} $\mu$. So, $T_{\lambda|\mu}$ implies
$T_{\lambda}$ {\em or} $T_{\mu}$ and $L^{\pm}_{\lambda|\mu}$ implies
$L^{\pm}_{\lambda}$ {\em or} $L^{\pm}_{\mu}$. If we define
$X^{+}_{\lambda}=C_{\lambda}$, $X^{+}_{\mu}=C_{\mu}$,
$X^{-}_{\lambda}=B_{\lambda}$, $X^{-}_{\mu}=B_{\mu}$, then the
spin-$\frac{1}{2}$ representation $\rho$ can be given as
\begin{eqnarray}
&&\rho (H_{\lambda})= \rho (H_{\mu})= 
\left( \begin{smallmatrix}1 & 0\\ 0 & -1\end{smallmatrix} \right)\\
&&\rho (X^{+}_{\lambda})= \rho (X^{+}_{\mu})= 
\left( \begin{smallmatrix}0 & 0\\ 1 & 0\end{smallmatrix} \right)\\
&&\rho (X^{-}_{\lambda})= \rho (X^{-}_{\mu})=
\left( \begin{smallmatrix}0 & 1\\ 0 & 0\end{smallmatrix} \right)\\
&&\rho (H'_{\lambda})= \rho (H'_{\mu})= 
\left( \begin{smallmatrix}1 & 0\\ 0 & 1\end{smallmatrix} \right)
\end{eqnarray}
and
\begin{equation}
\rho ({({X}^{\pm}_{\lambda|\mu})}^{2})=0 \qquad \qquad
\rho ({H}^{2}_{\lambda|\mu}) = 1 = \rho ({H'}^{2}_{\lambda|\mu})
\end{equation}
For vanishing colour parameters, the coloured $L^{\pm}$ functionals
reduce to the ordinary $L^{\pm}$ functionals for $\GLqtwo$.

\subsection{Coloured $RLL$ relations}

Similar to the uncoloured case, the commutation relations of the
algebra dual to a coloured quantum group can be obtained using the
modified or the {\em coloured} $RLL$ relations using the coloured
$L^{\pm}$ functionals of the previous section. Bearing close resemblence
to the coloured $RTT$ relations, these can be defined as
\begin{eqnarray}
R_{12} L^{\pm}_{2\lambda} L^{\pm}_{1\mu} &=& 
L^{\pm}_{1\mu} L^{\pm}_{2\lambda} R_{12}\\
R_{12} L^{+}_{2\lambda} L^{-}_{1\mu} &=&
L^{-}_{1\mu} L^{+}_{2\lambda} R_{12}
\end{eqnarray}
where $L^{\pm}_{1\mu}=L^{\pm}_{\mu}\ot \ident$ and
$L^{\pm}_{2\lambda}=\ident \ot L^{\pm}_{\lambda}$. Using the above
formulae, we obtain the commutation relations between the generating
elements of the algebra dual to the coloured quantum group
$GL_{q}^{\lambda,\mu}(2)$.

\begin{equation}
\begin{array}{lll}
\left[A_{\lambda},B_{\mu}\right]=B_{\mu} \qquad &
\left[D_{\lambda},B_{\mu}\right]=-B_{\mu} \qquad &\\
\left[A_{\lambda},C_{\mu}\right]=-C_{\mu} \qquad &
\left[D_{\lambda},C_{\mu}\right]=C_{\mu} \qquad &\\
\left[A_{\lambda},D_{\mu}\right]=0 \qquad &
\left[H_{\lambda},H_{\mu}\right]=0 \qquad &
\left[H'_{\lambda},\bullet\right]=0
\end{array}
\end{equation}
\begin{equation}
q^{-(\lambda + \mu)}C_{\lambda}B_{\mu} - 
q^{\lambda + \mu}B_{\mu}C_{\lambda} = 
\frac{q^{\lambda H_{\mu} + \mu H_{\lambda}}}{q-q^{-1}}
\left[
q^{-\frac{1}{2}(H_{\lambda}+H_{\mu})} 
q^{\lambda H'_{\lambda} - \mu H'_{\mu}} -
q^{\frac{1}{2}(H_{\lambda}+H_{\mu})}
q^{-\lambda H'_{\lambda} + \mu H'_{\mu}}
\right]
\end{equation}
\begin{equation}
\begin{array}{lll}
A_{\lambda}A_{\mu} &=& A_{\mu}A_{\lambda}\\
B_{\lambda}B_{\mu} &=& q^{2(\mu-\lambda)}B_{\mu}B_{\lambda}\\
C_{\lambda}C_{\mu} &=& q^{2(\lambda-\mu)}C_{\mu}C_{\lambda}\\
D_{\lambda}D_{\mu} &=& D_{\mu}D_{\lambda}
\end{array}
\end{equation}
where $H_{\lambda}$ and $H'_{\lambda}$ are as before. The relations
satisfy the $\lambda \leftrightarrow \mu$ exchange symmetry. The
associated coproduct of the elements of the dual algebra is given by
\begin{eqnarray}
\cop(A_{\lambda}) &=& A_{\lambda}\ot \ident + \ident\ot A_{\lambda}\\
\cop(B_{\lambda}) &=& B_{\lambda}\ot q^{A_{\lambda}-D_{\lambda}} +
\ident\ot B_{\lambda}\\
\cop(C_{\lambda}) &=& C_{\lambda}\ot q^{A_{\lambda}-D_{\lambda}} +
\ident\ot C_{\lambda}\\
\cop(D_{\lambda}) &=& D_{\lambda}\ot \ident + \ident\ot D_{\lambda}
\end{eqnarray}
The counit $\cnt(Y_{\lambda}) = 0;$ where $Y_{\lambda} = \{ A_{\lambda},
B_{\lambda}, C_{\lambda}, D_{\lambda} \}$ and the antipode is given as
\begin{eqnarray}
S(A_{\lambda}) &=& - A_{\lambda}\\
S(B_{\lambda}) &=& - B_{\lambda}q^{-(A_{\lambda}-D_{\lambda})}\\
S(C_{\lambda}) &=& - C_{\lambda}q^{-(A_{\lambda}-D_{\lambda})}\\
S(D_{\lambda}) &=& - D_{\lambda}
\end{eqnarray}
The Hopf algebra structure is also invariant under the $\lambda
\leftrightarrow \mu$ exchange symmetry. Thus we have defined a new
single-parameter coloured quantum algebra corresponding to $gl(2)$, which
in the monochromatic limit defines the standard uncoloured two-parameter
quantum algebra for $gl(2)$. Note that while the Hopf algebra structure
underlying a coloured quantum group and its dual is infinite dimensional,
it can be couched in a more familiar form of the finite dimensional case
as shown above. This makes it plausible to investigate the topological
aspect of duality for coloured quantum groups along the lines of the
treatment furnished in [15] for infinite dimensional Hopf algebras.

\section{Differential Calculus on $GL_{q}^{\lambda,\mu}(2)$}

We now proceed to the construction of a differential calculus on
the coloured
quantum group $GL_{q}^{\lambda,\mu}(2)$. Our preferred approach to this
problem is the $R$-matrix formalism, which we shall now generalise to the
case of coloured quantum groups.

\subsection{One-forms}

Analogous to the standard uncoloured quantum group, a bimodule $\Gamma$
(space of quantum one-forms $\omega$) is characterised by the commutation
relations between $\omega$ and $a_{\lambda(\mu)} \in\ca$, the {\em
coloured} quantum group or Hopf algebra under consideration,
\begin{equation}
\omega a_{\lambda(\mu)} = (\ident\ot f_{\lambda,\mu}) \cop
(a_{\lambda(\mu)}) \omega
\end{equation}
where $a_{\lambda(\mu)}$ denotes $a_{\lambda}$ ({\em respectively}
$a_{\mu}$) and the linear functional $f_{\lambda,\mu}$ is now defined in
terms of the coloured $L^{\pm}$ matrices as
\begin{equation}
f_{\lambda,\mu} = S(L_{\lambda|\mu}^{+})L_{\lambda|\mu}^{-}
\end{equation}
For the coloured quantum group $GL_{q}^{\lambda,\mu}(2)$, these read
\begin{equation}
S(L^{+}_{\lambda(\mu)}) = {(c^{+}{q}^{1/2})}^{-1} \begin{pmatrix}
q^{-H_{\lambda(\mu)}/2}q^{-\mu H_{\lambda(\mu)}+\lambda 
H'_{\lambda(\mu)}} & S(l^{+}_{\lambda(\mu)})_{12}\\
0 & q^{H_{\lambda(\mu)}/2}q^{-\mu H_{\lambda(\mu)}-\lambda
H'_{\lambda(\mu)}}
\end{pmatrix}   
\end{equation}
with
\begin{equation*}
S(l^{+}_{\lambda(\mu)})_{12}=
-q^{-1/2}(q-q^{-1}) q^{-H_{\lambda(\mu)}/2}q^{-\mu
H_{\lambda(\mu)}+\lambda H'_{\lambda(\mu)}} C_{\lambda(\mu)}
q^{H_{\lambda(\mu)}/2}q^{-\mu H_{\lambda(\mu)}-\lambda
H'_{\lambda(\mu)}}
\end{equation*}
\begin{equation}
L^{-}_{\lambda(\mu)} = c^{-} q^{-1/2} \begin{pmatrix}
q^{-H_{\lambda(\mu)}/2}q^{\lambda H_{\lambda(\mu)}-\mu H'_{\lambda(\mu)}}
& 0\\ q^{1/2}(q^{-1}-q)B_{\lambda(\mu)} & 
q^{H_{\lambda(\mu)}/2}q^{\lambda H_{\lambda(\mu)}+\mu H'_{\lambda(\mu)}}
\end{pmatrix}
\end{equation}
where $H_{\lambda(\mu)}=A_{\lambda(\mu)}-D_{\lambda(\mu)}$ and
$H'_{\lambda(\mu)}=A_{\lambda(\mu)}+D_{\lambda(\mu)}$
Thus, we have
\begin{equation}
\omega a_{\lambda(\mu)} = [(\ident\ot
S(L_{\lambda|\mu}^{+})L_{\lambda|\mu}^{-}) \cop (a_{\lambda(\mu)})]
\omega
\end{equation}
In terms of components, this can be written as
\begin{equation}
\omega_{ij} a_{\lambda(\mu)} = [(\ident\ot
S(l^{+}_{(\lambda|\mu)ki})l^{-}_{(\lambda|\mu)jl}) \cop
(a_{\lambda(\mu)})] \omega_{kl}
\end{equation}
using the expressions $L^{\pm} = l^{\pm}_{ij}$ and $\omega =
\omega_{ij}$ where $i, j = 1,2$. Using these relations, we obtain the
commutation relations of all the left-invariant one-forms with the
elements of the coloured quantum group $GL_{q}^{\lambda,\mu}(2)$ as
follows
\begin{eqnarray}
\omega^{1} a_{\lambda(\mu)} &=& \ms q^{-2+2(\lambda-\mu)}
a_{\lambda(\mu)}\omega^{1}\\
\omega^{+} a_{\lambda(\mu)} &=& \ms q^{-1+2\lambda} 
a_{\lambda(\mu)}\omega^{+}\\
\omega^{-} a_{\lambda(\mu)} &=& \ms q^{-1-2\mu}a_{\lambda(\mu)}
\omega^{-} + \ms (q^{-2}-1)q^{(\lambda-\mu)}b_{\lambda(\mu)}\omega^{1}\\
\omega^{2} a_{\lambda(\mu)} &=& \ms a_{\lambda(\mu)} \omega^{2} +
\ms (q^{-1}-q)q^{(\lambda+\mu)}b_{\lambda(\mu)}\omega^{+}
\end{eqnarray}
\begin{eqnarray}
\omega^{1} b_{\lambda(\mu)} &=& \ms b_{\lambda(\mu)} \omega^{1}\\
\omega^{+} b_{\lambda(\mu)} &=& \ms q^{-1+2\mu}b_{\lambda(\mu)}
\omega^{+} + \ms (q^{-2}-1)q^{(\lambda-\mu)}a_{\lambda(\mu)}\omega^{1}\\
\omega^{-} b_{\lambda(\mu)} &=&
\ms q^{-1-2\lambda}b_{\lambda(\mu)}\omega^{-}\\
\omega^{2} b_{\lambda(\mu)} &=& \ms q^{-2+2(\mu-\lambda)}
b_{\lambda(\mu)}\omega^{2} +
\ms (q^{-1}-q)q^{-(\lambda+\mu)}a_{\lambda(\mu)}\omega^{-}\\
& & + \ms (q^{-1}-q)^{2}b_{\lambda(\mu)}\omega^{1} \notag
\end{eqnarray}
\begin{eqnarray}
\omega^{1} c_{\lambda(\mu)} &=& \ms q^{-2+2(\lambda-\mu)}  
c_{\lambda(\mu)}\omega^{1}\\
\omega^{+} c_{\lambda(\mu)} &=& \ms q^{-1+2\lambda}
c_{\lambda(\mu)}\omega^{+}\\
\omega^{-} c_{\lambda(\mu)} &=& \ms q^{-1-2\mu}c_{\lambda(\mu)}
\omega^{-} + \ms (q^{-2}-1)q^{(\lambda-\mu)}d_{\lambda(\mu)}\omega^{1}\\
\omega^{2} c_{\lambda(\mu)} &=& \ms c_{\lambda(\mu)} \omega^{2} +
\ms (q^{-1}-q)q^{(\lambda+\mu)}d_{\lambda(\mu)}\omega^{+}
\end{eqnarray}
\begin{eqnarray}
\omega^{1} d_{\lambda(\mu)} &=& \ms d_{\lambda(\mu)} \omega^{1}\\
\omega^{+} d_{\lambda(\mu)} &=& \ms q^{-1+2\mu}d_{\lambda(\mu)}
\omega^{+} + \ms(q^{-2}-1)q^{(\lambda-\mu)}c_{\lambda(\mu)}\omega^{1}\\
\omega^{-} d_{\lambda(\mu)} &=&
\ms q^{-1-2\lambda}d_{\lambda(\mu)}\omega^{-}\\
\omega^{2} d_{\lambda(\mu)} &=& \ms q^{-2+2(\mu-\lambda)}
d_{\lambda(\mu)}\omega^{2} +
\ms(q^{-1}-q)q^{-(\lambda+\mu)}c_{\lambda(\mu)}\omega^{-}\\
& & + \ms (q^{-1}-q)^{2}d_{\lambda(\mu)}\omega^{1} \notag
\end{eqnarray}

where $\omega^{1}=\omega_{11}$, $\omega^{+}=\omega_{12}$,
$\omega^{-}=\omega_{21}$, $\omega^{2}=\omega_{22}$ and
$\ms=(c^{+})^{-1}c^{-}$.

\subsection{Vector fields}

The left-invariant vector fields $\chi_{ij}$ on $\ca$ are given by the
expression
\begin{equation}
\chi_{ij} = S(l^{+}_{(\lambda|\mu)ik})l^{-}_{(\lambda|\mu)kj} -
\delta_{ij}\cnt
\end{equation}
or simply
\begin{equation}
\chi = S(L^{+}_{\lambda|\mu})L^{-}_{\lambda|\mu} - \ident\cnt
\end{equation}
On elements of $GL_{q}^{\lambda,\mu}(2)$, the vector fields act as
\begin{eqnarray}
\chi_{ij} a_{\lambda(\mu)} &=&
(S(l^{+}_{(\lambda|\mu)ik})l^{-}_{(\lambda|\mu)kj} -
\delta_{ij}\cnt)a_{\lambda(\mu)}\\
\chi_{ij} a_{\lambda(\mu)} &=& \langle
S(l^{+}_{(\lambda|\mu)ik})l^{-}_{(\lambda|\mu)kj},
a_{\lambda(\mu)}\rangle - \delta_{ij}\cnt (a_{\lambda(\mu)})
\end{eqnarray}
for $a_{\lambda(\mu)}\in GL_{q}^{\lambda,\mu}(2)$. We obtain explicitly
the following

\begin{equation}
\begin{array}{ll}
\chi_{1}(a_{\lambda|\mu}) = \ms q^{-2+2(\lambda-\mu)}-1 & \qquad
\chi_{1}(b_{\lambda|\mu}) = 0\\
\chi_{+}(a_{\lambda|\mu}) = 0 & \qquad \chi_{+}(b_{\lambda|\mu}) = 0\\
\chi_{-}(a_{\lambda|\mu}) = 0 & \qquad \chi_{-}(b_{\lambda|\mu}) =
\ms (q^{-1}-q)q^{-(\lambda+\mu)}\\
\chi_{2}(a_{\lambda|\mu}) = \ms -1 & \qquad \chi_{2}(b_{\lambda|\mu}) = 0
\end{array}
\end{equation}
\begin{equation}
\begin{array}{ll}
\chi_{1}(c_{\lambda|\mu}) = 0& \qquad
\chi_{1}(d_{\lambda|\mu}) = \ms (q^{-1}-q)^{2} + \ms -1\\
\chi_{+}(c_{\lambda|\mu}) = \ms (q^{-1}-q)q^{-(\lambda+\mu)} & \qquad
\chi_{+}(d_{\lambda|\mu}) = 0\\
\chi_{-}(c_{\lambda|\mu}) = 0 & \qquad \chi_{-}(d_{\lambda|\mu}) = 0\\
\chi_{2}(c_{\lambda|\mu}) = 0 & \qquad \chi_{2}(d_{\lambda|\mu}) =
\ms q^{-2+2(\mu-\lambda)}-1
\end{array}   
\end{equation}

where $\chi_{0}=\chi_{11}$, $\chi_{+}=\chi_{12}$, $\chi_{-}=\chi_{21}$,
$\chi_{2}=\chi_{22}$, and $\ms =(c^{+})^{-1}c^{-}$. The left convolution
products are given as

\begin{equation}
\begin{array}{ll}
\chi_{1} \ast a_{\lambda(\mu)} =
(\ms q^{-2+2(\lambda-\mu)}-1)a_{\lambda(\mu)}
& \chi_{1} \ast b_{\lambda(\mu)} =
(\ms (q^{-1}-q)^{2}+ \ms -1)b_{\lambda(\mu)}\\
\chi_{+} \ast a_{\lambda(\mu)} =
(\ms (q^{-1}-q)q^{(\lambda+\mu)})b_{\lambda(\mu)}
& \chi_{+} \ast b_{\lambda(\mu)} = 0\\
\chi_{-} \ast a_{\lambda(\mu)} = 0
& \chi_{-} \ast b_{\lambda(\mu)} =
(\ms (q^{-1}-q)q^{-(\lambda+\mu)})a_{\lambda(\mu)}\\
\chi_{2} \ast a_{\lambda(\mu)} = (\ms -1)a_{\lambda(\mu)}
& \chi_{2} \ast b_{\lambda(\mu)} =
(\ms q^{-2+2(\mu-\lambda)}-1)b_{\lambda(\mu)}
\end{array}
\end{equation}   
\begin{equation}
\begin{array}{ll}
\chi_{1} \ast c_{\lambda(\mu)} =
(\ms q^{-2+2(\lambda-\mu)}-1)c_{\lambda(\mu)}
& \chi_{1} \ast d_{\lambda(\mu)} =
(\ms (q^{-1}-q)^{2}+ \ms -1)d_{\lambda(\mu)}\\
\chi_{+} \ast c_{\lambda(\mu)} =
(\ms (q^{-1}-q)q^{(\lambda+\mu)})d_{\lambda(\mu)}
& \chi_{+} \ast d_{\lambda(\mu)} = 0\\
\chi_{-} \ast c_{\lambda(\mu)} = 0
& \chi_{-} \ast d_{\lambda(\mu)} =
(\ms (q^{-1}-q)q^{-(\lambda+\mu)})c_{\lambda(\mu)}\\
\chi_{2} \ast c_{\lambda(\mu)} = (\ms -1)c_{\lambda(\mu)}
& \chi_{2} \ast d_{\lambda(\mu)} =
(\ms q^{-2+2(\mu-\lambda)}-1)d_{\lambda(\mu)}
\end{array}
\end{equation}

\subsection{Exterior derivatives}

Using the formula $\mathbf{d} a_{\lambda(\mu)} =\sum_{i}(\chi_{i} \ast
a_{\lambda(\mu)}) \omega^{i}$ for $a_{\lambda(\mu)} \in \ca$, we obtain
the action of the exterior derivative on the generating elements of
$GL_{q}^{\lambda,\mu}(2)$

\begin{eqnarray}
\mathbf{d} a_{\lambda(\mu)} &=&
(\ms q^{-2+2(\lambda-\mu)}-1)a_{\lambda(\mu)}\omega^{1} +
\ms (q^{-1}-q)q^{\lambda+\mu}b_{\lambda(\mu)}\omega^{+}\\
& &+ (\ms -1)a_{\lambda(\mu)}\omega^{2} \notag \\
\mathbf{d} b_{\lambda(\mu)} &=& 
(\ms (q^{-1}-q)^{2}+\ms -1)b_{\lambda(\mu)}\omega^{1} +
\ms (q^{-1}-q)q^{-(\lambda+\mu)}a_{\lambda(\mu)}\omega^{-}\\
& &+ (\ms q^{-2+2(\mu-\lambda)}-1)b_{\lambda(\mu)}\omega^{2} \notag \\
\mathbf{d} c_{\lambda(\mu)} &=&
(\ms q^{-2+2(\lambda-\mu)}-1)c_{\lambda(\mu)}\omega^{1} +
\ms (q^{-1}-q)q^{\lambda+\mu}d_{\lambda(\mu)}\omega^{+}\\
& &+ (\ms -1)c_{\lambda(\mu)}\omega^{2} \notag \\
\mathbf{d} d_{\lambda(\mu)} &=&
(\ms (q^{-1}-q)^{2}+\ms -1)d_{\lambda(\mu)}\omega^{1} +
\ms (q^{-1}-q)q^{-(\lambda+\mu)}c_{\lambda(\mu)}\omega^{-}\\ 
& &+ (\ms q^{-2+2(\mu-\lambda)}-1)d_{\lambda(\mu)}\omega^{2} \notag
\end{eqnarray}

$\mathbf{d}\ca$ generates $\Gamma$ as a left
$\ca$-module. This then defines a first order differential calclulus
$(\Gamma, \mathbf{d})$ on $GL_{q}^{\lambda,\mu}(2)$. Note that because
the colour parameters $\lambda$ and $\mu$ are continuously varying, the
differential calculus obtained is infinite dimensional. It can be checked
that the differential calculus on the uncoloured single-parameter quantum
group $\GLqtwo$ is recovered in the colourless limit
$\lambda = \mu = 0$. Furthermore, in the monochromatic limit $\lambda=\mu
\neq 0$, the differential calculus reduces to that of the uncoloured
two-parameter quantum group $\GLpqtwo$.

\section{Conclusions}

In this work, we have given a description of the
duality between a coloured quantum group and a coloured quantum algebra
employing the $R$-matrix approach. In particular, we have derived
explicitly the algebra dual to the coloured quantum group
$GL_{q}^{\lambda,\mu}(2)$. Following the duality, we have also
constructed a differential calculus on $GL_{q}^{\lambda,\mu}(2)$ by
providing a coloured generalisation of the $R$-matrix (or the {\em
constructive}) procedure. Both, the dual algebra as well as the
differential calculus for $GL_{q}^{\lambda,\mu}(2)$ reduce to that of the
ordinary $\GLqtwo$ in the colourless limit and to the uncoloured
two-parameter quantum group in the monochromatic limit. The results
obtained are general enough to be applicable to multi-parametric and
higher dimensional coloured quantum groups.
\par
It would be interesting to give dual basis for a coloured quantum group
and to investigate differential calculi on its quantum planes.

\section*{Acknowledgments}

It is a pleasure to thank Professor John Madore for discussions.

\section*{References}

\begin{enumerate}

\item
Kundu, A. and Basu-Mallick, B.: {\it J. Phys.} {\bf A27}
(1994) 3091; Basu-Mallick, B.: {\it Mod. Phys. Lett.} {\bf A9}
(1994) 2733; {\it Intl. J. Mod. Phys.} {\bf A10} (1995) 2851.

\item
Quesne, C.: {\it J. Math. Phys.} {\bf 38} (1997) 6018; {\it ibid} {\bf 39}
(1998) 1199.

\item
Parashar, P.: {\it Lett. Math. Phys.} {\bf 45} (1998) 105.

\item
Bazhanov, V.V. and Stroganov, Yu. G.: {\it Theor. Math. Phys.} {\bf 62}
(1985) 253.

\item
Murakami. J.: {\it Proc. Intl. Conf. of Euler Mathematical School: Quantum
Groups (Leningrad)}, Lecture Notes in Physics (Springer-Berlin), 1990,
p. 350.

\item
Ge, M.L. and Xue, K.: {\it J. Phys.} {\bf A24}
(1991) L895; {\it J. Phys.} {\bf A26} (1993) 281.

\item
Burdik, C. and Hellinger, P.: {\it J. Phys.} {\bf A25} (1992) L1023.

\item
Ohtsuki, T.: {\it J. Knot Theor. Its Rami.} {\bf 2} (1993) 211.

\item
Dabrowski, L. and Parashar, P.: {\it Lett. Math. Phys.} {\bf 38}
(1996) 331.

\item
Parashar, D. and McDermott, R.J.: {\it J. Math. Phys.} {\bf 41}
(2000) 2403; {\tt math.QA/9911244}.

\item
Faddeev, L.D., Reshetikhin, N. Yu. and Takhtajan,
L.A.: {\it Len. Math. J.} {\bf 1} (1990) 193.

\item
Dobrev, V.K.: {\it J. Math. Phys.} {\bf 33} (1992) 3419; Dobrev, V.K. and
Parashar, P.: {\it J. Phys.} {\bf A26} (1993) 6991.

\item
Jurco, B.: {\it Lett. Math. Phys.} {\bf 22} (1991) 177; Preprint CERN-TH
9417/94, 1994.

\item
Aschieri, P. and Castellani, L.: {\it Intl. J. Mod. Phys.} {\bf A8}
(1993) 1667.

\item
Bonneau, P., Flato, M. and Pinczon, G.: {\it Lett. Math. Phys.} {\bf 26}
(1992) 75; Bonneau, P., Flato, M., Gerstenhaber, M. and Pinczon, G.: {\it
Commun. Math. Phys.} {\bf 161} (1994) 125.

\end{enumerate}

\end{document}